\documentclass[12pt]{amsart}
\usepackage{amssymb}

\numberwithin{equation}{section}

\newcommand{\fullref}[1]{\ref{#1} on page~\pageref{#1}}

\newcommand{\ndash}{\nobreakdash-\hspace{0pt}}
\newcommand{\Ndash}{\nobreakdash--}

\newcommand{\dd}{{\mathrm{d}}}

\DeclareMathOperator{\pr}{pr}

\newcommand{\ev}{\mathrm{ev}}

\newcommand{\uOmega}{{\underline\Omega}}

\newcommand{\lo}{\mathrm{L}}

\newcommand{\pish}{\pi^\#}
\newcommand{\tTsM}{T^*M_{(\pi,\phi)}}

\newtheorem{Thm}{Theorem}[section]

\newtheorem{Lem}[Thm]{Lemma}

\newtheorem*{Thm*}{Theorem}

\theoremstyle{remark}

\newtheorem*{Ack}{Acknowledgment}

\theoremstyle{definition}

\newtheorem{lem}[Thm]{Lemma}
\newtheorem{cor}[Thm]{Corollary}
\newtheorem{pro}[Thm]{Proposition}

\newtheorem{defi}[Thm]{Definition}

\newcommand{\braket}[2]{\left\langle{\,{#1}\,,\,{#2}\,}\right\rangle}

\newcommand{\Lie}[2]{{\left[{\,{#1}\,,\,{#2}\,}\right]}}

\newcommand{\de}{\partial}

\def\gpd{\,\lower1pt\hbox{$\longrightarrow$}\hskip-.24in\raise2pt
               \hbox{$\longrightarrow$}\,}

\newcommand{\toto}{\rightrightarrows}
\newcommand{\gm}{\Gamma }
  
\newcommand{\ra}{\rangle}
\newcommand{\la}{\langle}  
\newcommand{\backl}{\mathbin{\vrule width1.5ex height.4pt\vrule height1.5ex}}
\newcommand{\per}{\backl}     
\newcommand{\be}{\begin{eqnarray*}}
\newcommand{\ee}{\end{eqnarray*}}  
\newcommand{\smalcirc}{\mbox{\tiny{$\circ $}}}

\let\Vec=\overrightarrow
\let\ceV=\overleftarrow
\let\Hat=\widehat

\newcommand\qq{\em}

\newcommand\lmp[1]{{\qq Lett.\ Math.\ Phys.\ \bf #1}}

\newcommand\anm[1]{{\qq Ann.\ Math.\ \bf #1}}

\newcommand\ptps[1]{{\qq Prog.\ Theor.\ Phys.\ Suppl.\ \bf #1}}
\newcommand{\tomega}{\tilde{\omega}}

\begin{document} 

\title{Integration of twisted Poisson structures}

\author[A.~S.~Cattaneo]{Alberto~S.~Cattaneo}
\address{Institut f\"ur Mathematik, Universit\"at Z\"urich--Irchel,  
Winterthurerstrasse 190, CH-8057 Z\"urich, Switzerland}  
\email{asc@math.unizh.ch}

\author[P.~Xu]{Ping Xu}
\address{Department of Mathematics, Penn State University,
University Park, Pa 16802, USA}
\email{ping@math.psu.edu}

\thanks{A.~S.~C. acknowledges partial support of SNF Grant No.~20-63821.00.
P. X. acknowledges partial support  of NSF Grant DMS00-72171.}

\begin{abstract}
Poisson manifolds may be regarded as the infinitesimal form of symplectic groupoids.
Twisted Poisson manifolds considered by \v Severa and Weinstein \cite{SW} 
are a natural generalization of the former which also arises in string theory. 
In this note
it is proved that twisted Poisson manifolds are in bijection with a
(possibly singular) twisted version of
symplectic groupoids.
\end{abstract}

\maketitle

\section{Introduction}\label{sec-intro}
Poisson manifolds may be regarded as the infinitesimal form of symplectic groupoids \cite{CDW},
i.e., Lie groupoids endowed with a multiplicative symplectic form. Up to singularities, Poisson
manifolds may be integrated to symplectic groupoids as described in \cite{CaFeOb} (conditions
under which integration with no singularities is possible are given in \cite{CFPoiss}).
In this paper we generalize this result to the case when the two structures (of symplectic
groupoid and of Poisson manifold) are twisted by a closed $3$\ndash form.

Let $M$ be a smooth manifold. A pair $(\pi,\phi)$, where
$\pi$ is a bivector field and $\phi$ is a closed $3$\ndash form, is called
a \textsf{twisted Poisson structure} if it satisfies the equation
\begin{equation}\label{phipoisson}
\Lie\pi\pi=\frac12\wedge^3\pi^\#\phi,
\end{equation}
where $\Lie{\ }{\ }$ denotes the Schouten--Nijenhuis bracket and
$\pi^\#$ is the vector bundle homomorphism $T^*M\to TM$ induced by $\pi$
(viz., $\pi^\#(x)(\sigma):=\pi(x)(\sigma,\bullet)$, with $x\in M$,
$\sigma\in T^*_xM$). According to \cite{SW}, one also says that
$\pi$ is a \textsf{$\phi$-Poisson tensor}. 
In the case $\phi=0$ one recovers the usual notions of Poisson
tensor and Poisson manifold. Twisted Poisson structures have been
extensively studied in the physics literature, e.g., \cite{P, CS, KS}.

As explained in \cite{SW},
a twisted Poisson structure induces a Lie algebroid structure on
$T^*M$ with anchor map $\pi^\#$ and Lie bracket of sections
$\sigma$ and $\tau$ defined by
\begin{equation}\label{tLie}
\Lie\sigma\tau :=
\lo_{\pish\sigma}\tau-\lo_{\pish\tau}\sigma
-\dd\pi(\sigma,\tau)+\phi(\pish\sigma,\pish\tau,\bullet).
\end{equation}
In particular, $\forall f,g\in C^\infty(M)$ we have:
\begin{align}
\label{eq:twist1}
[\dd f, \dd g]&=\dd\{f, g\}+\phi (X_{f}, X_{g}, \bullet )\\
\intertext{and}
\label{eq:twist2}
[X_f, X_g]&=X_{\{f, g\} }+\pi^{\#}(\phi (X_{f}, X_{g}, \bullet ) ),
\end{align}
with $\{f,g\}=\pi(\dd f,\dd g)$ and
$X_f=\pi^\#\dd f$.

We will denote this Lie algebroid by $\tTsM$.
Sections of its exterior algebra are ordinary differential forms.
One may define a derivation $\delta$ deforming the de~Rham differential
by $\phi$; viz., we define a graded derivation
$\delta\colon \Omega^*(M)\to\Omega^{*+1}(M)$ by setting
$\delta f=\dd f$ if $f\in C^\infty(M)$ and
\[
\delta\sigma = \dd\sigma -\iota_{\pish\sigma}\phi,
\]
if $\sigma\in\Omega^1(M)$. It turns out that
\[
\delta\Lie\sigma\tau=\Lie{\delta\sigma}\tau+\Lie\sigma{\delta\tau},
\qquad\forall\sigma,\tau\in\Omega^1(M),
\]
and that $\delta^2=\Lie\phi{\bullet}$
(where we have extended the Lie bracket to the whole of $\Omega^*(M)$
as a biderivation). So $(\tTsM,\delta)$ constitutes an example
of a quasi Lie bialgebroid \cite{R, ILX}, a generalization of
Drinfeld's quasi Lie bialgebras \cite{dr:quasi, YKS}.

If $\tTsM$ may be integrated to a Lie groupoid $(G\toto M,\alpha,\beta)$
(i.e., if it exists a Lie groupoid $G$ whose Lie algebroid is $\tTsM$), 
the differential $\delta$ induces extra structure on $G$.
Namely, denoting by $\alpha$ and $\beta$ the source and target maps
of $G$,
then $G$ may be endowed with a non-degenerate, multiplicative $2$\ndash form
$\omega$ that satisfies
\[
\dd\omega=\alpha^*\phi-\beta^*\phi.
\]
In other words, $(\omega,\phi)$ is a $3$\ndash cocycle for the double
complex $\Omega^*(G^{(*)})$, where $G^{(0)}=M$, $G^{(1)}=G$ and elements of
$G^{(k)}$ are $k$\ndash tuples of elements of $G$ that may be multiplied
(in the given order). One differential is de~Rham and the other is the
groupoid-complex differential. Observe that in the true Poisson case
(i.e., $\phi=0$), $\omega$ is closed, so $G$ is an ordinary symplectic groupoid.
In the general
case, $G$ is called a \textsf{non-degenerate 
twisted symplectic groupoid},
and the non-degenerate $2$\ndash form $\omega$ is said to be
 \textsf{relatively $\phi$\ndash closed}. The main theorem of the paper 
(conjectured in \cite{SW}) is
\begin{Thm*}
There is a bijection between  integrable twisted Poisson structures 
and  source-simply connected non-degenerate
twisted symplectic groupoids.
\end{Thm*}
Here   ``integrable twisted Poisson structure" means 
that the associated Lie algebroid is integrable.

In Sect.~\ref{sec-qsg} we give an introduction 
to non-degenerate twisted symplectic groupoids
and prove that they induce 
 twisted Poisson structures on the base manifolds (Theorem~\fullref{thm-qsgtP}).

In Sect.~\ref{sec-proof}
we prove the Theorem, though in a more general setting.
In fact, as shown in the generaliztion \cite{CF} (see also \cite{S}) of the construction
given in \cite{CaFeOb}),
to any Lie algebroid $A$ one can
associate a topological source-simply-connected
groupoid $G(A)$, which is the Lie groupoid integrating $A$ whenever
$A$ is integrable. The topological groupoid $G(A)$ is defined as the leaf
space of a smooth foliation, as we recall in Sect.~\ref{sec-int};
so it makes sense to define on it a notion of smooth functions and forms.
In the case when $A$ is $\tTsM$, we prove that $G(A)$ may always
be endowed with a non-degenerate, multiplicative, 
relatively $\phi$\ndash closed $2$\ndash form $\omega$. 
The construction is a modification, described in Sect.~\ref{sec-qsr},
of the method developed in \cite{CaFeOb}, where 
the true Poisson case (i.e., $\phi=0$) was dealt with.

As a final remark, we mention that
 general multiplicative $2$\ndash forms, 
their infinitesimal counterparts and their integrations are  being treated
in \cite{BCWZ}.

\begin{Ack}
We  thank Jim Stasheff and Alan Weinstein for useful discussions and comments.
We   acknowledge  Zurich University (P. X.) and
 Penn State  University (A. S. C.) for their kind  hospitality during
the preparation of the work.  
\end{Ack}

\newcommand{\Palg}{P(\tTsM)}
\newcommand{\Galg}{G(\tTsM)}
\newcommand{\PGamma}{P_0\Gamma(\tTsM)}

\section{Non-degenerate twisted symplectic groupoids}\label{sec-qsg}

\begin{defi}
A  non-degenerate twisted symplectic groupoid is a
  Lie groupoid $({G}\toto {M}, \alpha , \beta )$
equipped with a non-degenerate $2$\ndash form 
$\omega \in \Omega^2 (G )$ and
a $3$\ndash form $\phi \in \Omega^3 (M)$ such that
\begin{enumerate}
\item $\dd \phi =0$;
\item $\dd \omega =\alpha^* \phi  -\beta^* \phi $;
\item $\omega $ is multiplicative, i.e., 
the $2$\ndash form $(\omega , \omega , -\omega )$
vanishes when  being restricted to  the graph of the
groupoid multiplication $\Lambda\subset   G\times G\times G$.
\end{enumerate} 
\end{defi}

Let $\pi_G$ denote the bivector field on $G$ corresponding
to $\omega$. Then $(\pi_G, \Omega )$, where
 $ \Omega=\alpha^* \phi  -\beta^* \phi $, defines a twisted
Poisson structure on $G$ in the sense of \cite{SW}.

 
For any $\xi \in \gm (A)$, by $\Vec{\xi}$ and $\ceV{\xi}$
we denote its corresponding right and left invariant 
vector fields on the groupoid $G$, respectively.
The following properties can be easily verified.
                                                       
\begin{pro}
\label{pro:2.2}
~
\begin{enumerate}
\item $\epsilon^* \omega= 0$, where $\epsilon\colon M\to G$ is
the natural embedding;
\item $i^* \omega= -\omega$, where $i\colon G \to G$ is the  groupoid
inversion;
\item  for any $\xi, \eta \in \gm (A)$, $\omega (\Vec{\xi} , \Vec{\eta})$
is  a right invariant function on $G$,
 and $\omega (\ceV{\xi} ,\ceV{\eta})$ is a left invariant function on $G$;
\item   $\omega (\Vec{\xi} ,\ceV{\eta})=0$;
\item $\omega (\Vec{\xi} , \Vec{\eta})(x) 
=-\omega (\ceV{\xi} ,\ceV{\eta})(x^{-1})$.
\end{enumerate}
\end{pro}
\begin{proof} The proof is standard, and essentially follows from
the multiplicativity of $\omega$.

\begin{enumerate}
\item For any $\delta'_{m}, \ \delta''_{m}\in T_m M$,
since $(\delta'_{m} , \delta'_{m} , \delta'_{m}),\ 
(\delta''_{m} , \delta''_{m}, \delta''_{m}) \in T\Lambda $,
it follows that $\omega (\delta'_{m} , \delta''_{m} )=0$.
\item $\forall x\in G$ and $\forall \delta'_{x}, \ \delta''_{x}
\in T_{x}G$, it is clear that
 $ (\delta'_{x} , i_{*} \delta'_{x} , \alpha_* \delta'_{x} ),
\ (\delta''_{x} , i_{*} \delta''_{x} , \alpha_* \delta''_{x} )\in 
T\Lambda $. Thus, by (1), we have
$$\omega (\delta'_{x}, \delta''_{x})+\omega (i_{*} \delta'_{x},
i_{*} \delta''_{x} )=0, $$
and (2) follows.
\item For any $\xi, \eta \in \gm (A)$,
 $( \Vec{\xi}(x), 0_{y}, \Vec{\xi}(xy)), \
( \Vec{\eta}(x), 0_{y}, \Vec{\eta }(xy)) \in T\Lambda$.
Thus
$$\omega (\Vec{\xi}(x),  \Vec{\eta}(x))-\omega (\Vec{\xi}(xy) , 
\Vec{\eta }(xy) )=0. $$
Hence $\omega (\Vec{\xi} , \Vec{\eta})$
is  a right invariant function on $G$.
  Similarly,  $\omega (\ceV{\xi} ,\ceV{\eta})$ is a left invariant function
 on $G$.  
\item By considering the vectors
$( \Vec{\xi}(x), 0_{\beta (x)}, \Vec{\xi}(x))$  and
$(0_x , \ceV{\eta}(\beta (x) ), \ceV{\eta}(x ))\in
T\Lambda $, we obtain
$\omega (\Vec{\xi}(x), \ceV{\eta}(x ))=0$.
\item follows from (2) and the fact that $i_* \Vec{\xi}=-\ceV{\xi}$.
\end{enumerate}
\end{proof}

Define a section $\gamma\in  \gm (\wedge^2 A^* )$ and
a bundle map: $\lambda\colon A \to T^*M$  by

\begin{equation}
\label{omegavv}
\omega (\Vec{\xi},  \Vec{\eta} )= \alpha^ *\gamma (\xi, \eta  ),
\ \ \forall \xi, \eta \in \gm (A),
\end{equation}   
and
\begin{equation}
\label{omegavh}
\la \lambda (\xi ), v\ra =\omega (\Vec{\xi} (m), v), \ \ \forall
\xi \in A|_m ,  v\in T_m M
\end{equation}

\begin{lem}
~
\begin{enumerate}
\item $ \omega (\ceV{\xi} , \ceV{\eta } )= -\beta^*  \gamma (\xi, \eta  ),
\ \ \forall \xi , \ \eta \in \gm (A)$;
\item for all $\xi , \eta \in \gm (A)$,
\begin{equation}
\label{eq:ga}
\gamma (\xi , \eta )=\la \rho (\xi), \lambda (\eta )\ra ;
\end{equation}
\item $\lambda\colon A \to T^*M$ is a  vector bundle isomorphism.
\end{enumerate} 
\end{lem}

\pagebreak
\begin{proof}
~
\begin{enumerate}
\item follows from Proposition~\ref{pro:2.2} (5).
\item We have
 $$\omega (\Vec{\xi},  \Vec{\eta} ) =\omega (\Vec{\xi} -\ceV{\xi},
\Vec{\eta})=\omega (\Vec{\eta}, \rho (\xi ))=\la \rho (\xi), \lambda (\eta )\ra $$
\item Assume that $ \lambda  (\xi )=0$. That is,
$\omega (\Vec{\xi} (m), v)=0, \ \forall v\in T_m M$, which implies that
$\Vec{\xi} (m) \per \omega =0$ by Proposition~\ref{pro:2.2} (4). 
Hence $\xi=0$ since
$\omega$ is non-degenerate.  This means that $ \lambda$ is injective.
On the other hand, assume that $v\in (\lambda (A|_m) )^{\perp}$.
Then $\omega (\Vec{\xi} (m), v)=0, \ \forall \xi \in A|_m$.
Thus $v\per  \omega =0$ using Proposition~\ref{pro:2.2} (1),
 which implies that
$v=0$. Therefore $\lambda $ is  surjective. 
\end{enumerate} \end{proof}

\begin{lem}
For any $f\in C^{\infty}(M)$,
\begin{equation}
\Vec{\lambda^{-1}(\dd f)}=X_{\alpha^* f}; \ \ \ 
\ceV{\lambda^{-1}(\dd f)}=X_{\beta^* f}.
\end{equation} 
\end{lem}
\begin{proof} First, one shows that $X_{\alpha^* f}$ is a right invariant
vector field  on $G$ and $X_{\beta^* f}$ is a left invariant 
vector field. This can be shown using the  same argument
as in the case of symplectic groupoids \cite{CDW}. Namely
the multiplicativity of $\omega$ together with dimension counting
implies that the graph $\Lambda$ is coisotropic with respect
to $(\pi_{G}, \pi_{G} , -\pi_{G})$. The later
implies that $X_{\alpha^* f}$ is a right invariant
vector field  on $G$ and $X_{\beta^* f}$ is a left invariant
vector field. 

Second, for any $v\in T_m M$, we have
$$\omega (X_{\alpha^* f}(m) , v)=\la \alpha^* \dd f(m), v\ra
=\la \dd f(m),   \alpha_* v\ra 
=\la \dd f(m),  v\ra$$
It thus follows that $\lambda (X_{\alpha^* f} )= \dd f$, or 
$\Vec{\lambda^{-1}(\dd f)}=X_{\alpha^* f}$.
The other equation can be proved similarly. \end{proof}

By pulling back the $2$\ndash form $\gamma \in \gm (\wedge^2 A^* )$ via
$\lambda^{-1}$, one obtains a bivector field $\pi \in \gm (\wedge^2 TM)$.
We introduce a bracket and Hamiltonian vector
fields by the usual definitions; i.e., 
$\{f ,g\}=\pi (\dd f, \dd g)$ and $X_{f}=\pi^{\#} (\dd f)$. 

\begin{cor}
\begin{equation}
\alpha_* \pi_G =\pi ; \ \ \ \beta_* \pi_G =-\pi;
\end{equation} 
or equivalently
\begin{equation} 
\alpha_* X_{\alpha^* f}=X_f ; \ \ \ \beta_* X_{\beta^*f}=-X_f, \ \ \forall
f\in C^{\infty}(M).
\end{equation}
\end{cor}   
\begin{proof}  For any $f , \ g \in C^{\infty}(M)$, 
\begin{multline*}
\{  \alpha^* f , \alpha^* g\}=
\omega (X_{\alpha^* f}, X_{\alpha^* g} )
= \omega (\Vec{\lambda^{-1}(\dd f)}, \Vec{\lambda^{-1}(\dd g)} )=\\
=\alpha^* (\pi (\dd f, \dd g ))=\alpha^* \{f, g\}.
\end{multline*}
Similarly, we have $\{  \beta^* f ,  \beta^* g\}=
-\beta^* \{f, g\}$.
\end{proof}   

We are now ready to prove the main result of the section.

\begin{Thm}\label{thm-qsgtP}
~
\begin{enumerate}
\item $\pi$ is a $\phi$-Poisson tensor in the sense of 
\cite{SW}---i.e., it satisfies \eqref{phipoisson}.
\item The bundle map $\lambda\colon A\to T^*M$ establishes a Lie algebroid isomorphism, where
the Lie algebroid on $T^*M$ is induced by the twisted Poisson
tensor $\pi$ as given by Eq.\  (\ref{tLie}).
\end{enumerate}
\end{Thm}
\begin{proof} Let $\Omega =\alpha^* \phi -\beta^* \phi$.
Thus $\forall f, \ g\in C^{\infty}(M)$
\be
(X_{\alpha^* f}\wedge X_{\alpha^* g}) \per \Omega 
&=&(X_{\alpha^* f}\wedge X_{\alpha^* g}) \per\alpha^* \phi\\
&=&\alpha^* [(\alpha_* X_{\alpha^* f} \wedge \alpha_* X_{\alpha^* g})
\per \phi ]\\
&=&\alpha^* [ X_{f}\wedge X_g \per \phi ].
\ee

Thus by Eq.\  (\ref{eq:twist2}) 

\be
[X_{\alpha^* f}, X_{\alpha^* g}]-X_{\{\alpha^* f , \alpha^* g\}}
=\pi^{\#}_{G } ( \Omega (X_{\alpha^* f}, X_{\alpha^* g}, \bullet )\\
=\pi^{\#}_{G } (\alpha^* \phi ( X_{f}, X_g , \bullet ))
\ee

Thus it follows that
$$\lambda [X_{\alpha^* f}, X_{\alpha^* g}]
=\dd\{f, g\}+\phi ( X_{f}, X_g , \bullet ).$$
Note that $\lambda$ intertwines  the anchors: 
$\pi^{\#}\smalcirc \lambda =\rho$, according to Eq.\  (\ref{eq:ga}).
Therefore, using Lie algebroid properties, one
shows that the push forward Lie algebroid
on $T^*M$ via $\lambda$ is given by Eq.\  (\ref{tLie}). 
This forces, by the Jacobi identity, $\pi$ to be
$\phi$-Poisson, and $\lambda $ is a Lie algebroid isomorphism
between $A$ and $(T^*M)_{\pi , \phi }$. \end{proof}

\section{Integration of $\tTsM$}\label{sec-int}
We briefly describe the integration procedure for Lie algebroids
of \cite{CF,S}, adapted to the case of $\tTsM$. First one defines
the manifold $P(\tTsM)$ of $C^1$-Lie algebroid morphisms $TI\to\tTsM$,
where
$I$ is the interval $[0,1]$ and $TI$ is given its canonical Lie algebroid
structure. An element of $P\tTsM$ consists of a $C^2$-path $X\colon I\to M$
together with a section $\eta$ of $T^*I\otimes X^*T^*M$ 
satisfying
\[
\dd X = \pish(X)\eta.
\]
On this manifold one may consider as equivalent two elements which are
related by a Lie algebroid morphism $T(I\times I)\to \tTsM$ that fixes
the endpoints. The quotient space $G(\tTsM)$ may be given a groupoid
structure. For our purposes it is however better to use a different
description of $G(\tTsM)$, i.e., as the leaf space of a foliation.
Namely, let $P_0\Gamma(\tTsM)$ be the space of $C^2$-paths
in the Lie algebra of sections of $\tTsM$ with endpoints at zero.
We give this space the structure of a Lie algebra by the pointwise
Lie bracket. One may then define an infinitesimal action of this Lie
algebra on $P(\tTsM)$. To describe it, we prefer to introduce local
coordinates $\{x^i\}$
on $M$ (alternatively, one may use a torsion-free connection).
Since $\{\dd x^i\}$ is a local basis of sections of $\tTsM$,
we may define structure functions $f$ by
\[
\Lie{\dd x^i}{\dd x^j}=f^{ij}_k\,\dd x^k,
\]
where a sum over repeated indices is understood. If we write locally
$\pi=\pi^{ij}\de_i\de_j$ and $\phi=\phi_{ijk}\dd x^i\dd x^j\dd x^k$,
we may compute:
\[
f^{ij}_k=\de_k\pi^{ij} + \pi^{mi}\pi^{nj}\phi_{mnk}.
\]
The action is then as follows. To $B\in P_0\Gamma(\tTsM)$ we associate
a vector field $\xi_B$ on $\Palg$. We can always write
$\xi_B=\xi_B^h+\xi_B^v$ with
$\xi_B^h(X,\eta)\in\Gamma(I,X^*TM)$ and
$\xi_B^v(X,\eta)\in\Gamma(I,T^*I\otimes X^*T^*M)$. We set then
\begin{subequations}\label{xib}
\begin{align}
(\xi_B^h(X,\eta))^i &= -\pi^{ij}(X)\,(B_X)_j,\\
(\xi_B^v(X,\eta))_i &= -\dd (B_X)_i - f^{rs}_i(X)\,\eta_r\,(B_X)_s,
\end{align}
\end{subequations}
where $B_X$ is the section of $X^*T^*M$ defined by $B_X(t)=B(t)(X(t))$.

Thus, the infinitesimal action of $\PGamma$ defines a foliation
on $\Palg$ and $\Galg$ is its quotient space. Let us briefly recall
its groupoid structure. The target map $\alpha$ associates to
a class of morphisms $(X,\eta)$ the value of $X$ at $0$, while
the source map $\beta$ associates to it the value of $X$ at $1$
(observe that the infinitesimal action preserves the endpoints of $X$).
The identity section associates to a point $m$ in $M$ the class $\epsilon(m)$
of the constant path at $m$ with $\eta=0$. The product is obtained
by joining the base paths and restricting the fiber maps consequently
(the product is more precisely defined on smooth
representatives such that $\eta$ vanishes with its derivatives at the
endpoints).

\section{Quasi-symplectic reduction}\label{sec-qsr}
In this section we describe how to obtain $\Galg$ by some
sort of symplectic reduction, though our replacement for a symplectic
form will be a non-degenerate but not necessarily closed $2$\ndash form.

Let $T^*PM$ denote the manifold of $C^1$-bundle maps $TI\to T^*M$
(over $C^2$\ndash maps). This space is morally a cotangent bundle
and as such it has a canonical symplectic structure $\Omega_0$.
Explicitly, a point in $T^*PM$ is a pair $(X,\eta)$ where $X$ is a
$C^2$-path $I\to M$ and $\eta$ is 
a $C^1$-section  of $T^*I\otimes X^*T^*M$. The tangent space
at $(X,\eta)$ is the direct sum of
$T^h_{(X,\eta)}T^*PM=\Gamma(I,X^*TM)$
and
$T^v_{(X,\eta)}T^*PM=\Gamma(I,T^*I\otimes X^*T^*M)$.
Using this splitting, we write
\begin{equation}
\label{Omega0}
\Omega_0(X,\eta)(\xi_1\oplus e_1,\xi_2\oplus e_2)=
\int_I \braket{e_1}{\xi_2} - \braket{e_2}{\xi_1},
\end{equation}
where $\braket{\ }{\ }$ denotes the canonical pairing between
tangent and cotangent fibers of $M$.

Using the $3$\ndash form $\phi$ on $M$ we may also define
a second $2$\ndash form on $T^*PM$:
\begin{equation}
\label{Omega1}
\Omega_1(X,\eta)(\xi_1\oplus e_1,\xi_2\oplus e_2)=
\frac12 \int_I \phi(X)(\pish(X)\eta,\xi_1,\xi_2).
\end{equation}
The $2$\ndash form $\Omega=\Omega_0+\Omega_1$ is still non-degenerate but
no longer closed.

The manifold $\Palg$ introduced in the previous Section may be regarded
as a submanifold of $T^*PM$. If we introduce ``momentum maps''
$H\colon T^*PM\to \PGamma^*$ by
\[
H_B(X,\eta)=\int_I\braket{B_X}{\dd X-\pish(X)\eta},
\]
then $\Palg$ is $H^{-1}(0)$. One may check that $\dd H_B$ lies
in the image of $\Omega$ for any $B\in\PGamma$; so, since $\Omega$
is non-degenerate, one
may define a map $B\to\hat\xi_B$ that associates a vector field
$\hat\xi_B$ on $T^*PM$ to $B$ by
\begin{equation}\label{ham}
\iota_{\hat\xi_B}\Omega = \dd H_B.
\end{equation}
One may easily
check that the restriction of $\hat\xi_B$ to $\Palg$ is tangent to it.
More to the point, one may check that the vector field
on $\Palg$ so obtained is precisely the $\xi_B$ of \eqref{xib}
which defines the infinitesimal action of $\PGamma$ on $\Palg$.


\section{Proof of the Theorem}\label{sec-proof}
In the setting of the previous Section, we want to prove that
the restriction $\uOmega$
of $\Omega$ to $\Palg$ is basic w.r.t.\ to the
projection $p\colon\Palg\to \Galg$, viz., $\uOmega=p^*\omega$; moreover,
we want to prove that $\omega$ satisfies all the required conditions.

Observe that $\uOmega$ is automatically
horizontal by \eqref{ham}. On the other hand, unlike the usual symplectic
case, it is not clear that $\uOmega$ is also invariant; in fact,
at first, we may only see that
$\lo_{\xi_B}\uOmega=\iota_{\xi_B}\dd\uOmega=\iota_{\xi_B}\dd\uOmega_1$,
where $\uOmega_1$ denotes the restriction of $\Omega_1$ to $\Palg$.
To proceed, we must understand $\uOmega_1$ better.

Let $PM$ be the manifold of $C^2$-paths in $M$. Let
$\ev\colon I\times PM\to M$ be the evaluation map and
$\pr\colon I\times PM\to PM$ the projection to the second factor.
Define $\Phi=\pr_*\ev^*\phi\in\Omega^2(PM)$, where $\pr_*$ denotes
integration along the fiber. If we finally denote by
$q\colon \Palg\to PM$ the map that retains only the base map
of the Lie algebroid morphism, we realize immediately that
\[
\uOmega_1 = q^*\Phi.
\]
By the generalized Stokes' Theorem and the fact that $\phi$ is closed,
we obtain $\dd\Phi=\underline\alpha^*\phi-\underline\beta^*\phi$,
where $\underline\alpha$ and $\underline\beta$ are the maps $PM\to M$
that assign to a path its values at $0$ and at $1$, respectively.
Thus,
\[
\dd\uOmega=q^*(\underline\alpha^*\phi-\underline\beta^*\phi).
\]
Since the vector field $\xi_B$ does not move the endpoints, we conclude
that $\iota_{\xi_B}\dd\uOmega=0$, viz., that $\uOmega$ is invariant as well.
We write then $\uOmega=p^*\omega$ as at the beginning of the Section.
The $2$\ndash form $\omega$ on $\Galg$ is clearly multiplicative since
the product is defined by joining the paths and $\Omega$ is defined
as an integral. Moreover, recalling the definition of the source
and target map $\beta$ and $\alpha$, we observe that
$\underline\alpha\circ q=\alpha \circ p$ and
$\underline\beta\circ q=\beta\circ p$.
So we may write the equation above as
\[
\dd\uOmega=p^*(\alpha^*\phi-\beta^*\phi).
\]
Since $\dd\uOmega=p^*\dd\omega$ and $p$ is a surjection, 
this shows that $\omega$ is relatively $\phi$\ndash closed.


Finally, we need to prove that the $2$\ndash form $\omega$ is non-degenerate.
It is clear from the construction that $\omega$ is non-degenerate
along the identity $M$. The claim thus follows from the following:

\begin{Lem}
A multiplicative $2$\ndash form $\omega \in \Omega^2 (G )$ on
a Lie groupoid $G\toto M$ is non-degenerate if and only if
it is non-degenerate along the identity $M$.
\end{Lem}
\begin{proof} First of all,  note that for any $\delta_x \in T_x G$, and
$\xi\in \gm (A)$, we have
\begin{eqnarray}
&&\omega (\ceV{\xi}(x), \delta_x )=\omega  (\ceV{\xi}(v), \beta_{*}\delta_x ) \label{eq:1}\\ 
&&\omega (\Vec{\xi}(x), \delta_x )=\omega  (\Vec{\xi}(u), \alpha_{*}\delta_x ),  \label{eq:2}
\end{eqnarray}
where  $u=\alpha (x)$ and $v=\beta (x)$. Eq.\  (\ref{eq:1}), for
instance,  follows from the
fact that both $(\delta_x ,\delta_x , \beta_{*}\delta_x )$,  and
$(0, \ceV{\xi}(x), \ceV{\xi}(v))$ are tangent to the graph of  the groupoid
multiplication $\Lambda \subset G\times G \times G$.
Eq.\  (\ref{eq:2}) can be proved similarly. 
Now assume that $\delta_x \in \ker \omega_x$. It follows from 
Eq.\  (\ref{eq:1}) that $\beta_{*}\delta_x \in \ker \omega_v$ since $M$ is
isotropic with respect to $\omega$. Therefore $\beta_{*}\delta_x =0$ by assumption.
Hence $\delta_x =\Vec{\eta}(x) $. On the other hand, according
to Eq.\  (\ref{eq:2}),  one has $\omega  (\Vec{\eta}(u), T_u M )=0$ since
$\alpha$ is a submersion. This implies that $\Vec{\eta}(u)\in \ker \omega_u$.
Therefore $\Vec{\eta}(u) =0$ by assumption. This implies that $\delta_x =\Vec{\eta}(x) 
=0$. This concludes the proof. \end{proof}

We need now  to prove that  the correspondence between
$\phi$-twisted Poisson structures and   twisted symplectic
groupoids is a bijection. The proof is divided into two steps.

\subsection*{Step 1}
By construction (see \cite{CaFeOb,CF}) the Lie algebroid of $G(\tTsM)$ is
$\tTsM$.
As discussed in Sect.~\ref{sec-qsg},
the 
relatively $\phi$\ndash closed, multiplicative, non-degenerate $2$\ndash form 
$\omega $ determines
an automorphism $\lambda$ of $T^*M$ and a bivector field $\gamma$ on $M$ as in Eqs.\ \eqref{omegavv}
and \eqref{omegavh}. We have to show that $\lambda$ is the identity 
and that
$\gamma=\pi$. First of all we observe that it is enough to consider \eqref{omegavv} at
the unit element $\epsilon(m)\in G(\tTsM)$ corresponding to $m\in M$:
 \[
\omega(\epsilon(m))(\Vec{\xi_1}(\epsilon(m)),  \Vec{\xi_2}(\epsilon(m)) )= \gamma(m) (\xi
_1, \xi_2),
\ \ \forall \xi_1, \xi_2 \in A|_m.
\]
By construction $\epsilon(m)$
is the equivalence class of the  path $X(t)=m$, $\eta(t)=0$,
$\forall t\in I=[0,1]$. The vector field $\Vec{\xi_i}$, $i=1,2$, evaluated at $\epsilon(m
)$ is
the projection to $T_{\epsilon(m)}G(\tTsM)$ of the vector $\Hat\xi_i\in T_{(m,0)}P(\tTsM)
$ defined by
$\Hat\xi_i(t)=(\pi^\#(m)\xi_i\,t,\xi_i\dd t)$.
 Observing then that for $\eta=0$ the $2$\ndash form
$\Omega_1$ of Eq.\ \eqref{Omega1} vanishes, we get, also using \eqref{Omega0},
\begin{multline*}
\omega(\epsilon(m))(\Vec{\xi_1}(\epsilon(m)),  \Vec{\xi_2}(\epsilon(m)) )=
\Omega_0(m,0)(\Hat\xi_1,\Hat\xi_2) =\\
=2\int_0^1 \pi(m)(\xi_1,\xi_2)\,t\;\dd t =  \pi(m)(\xi_1,\xi_2),
\end{multline*}
which shows $\gamma=\pi$. As for \eqref{omegavh}, observe that $\omega(\epsilon(m))(\Vec{
\xi_1}(\epsilon(m)),v)$ is
just $\Omega_0(m,0)(\Hat\xi_1,\Hat v)$ with $\Hat v(t)=(v,0)$. As a consequence,
\[
\omega(\epsilon(m))(\Vec{\xi_1}(\epsilon(m)),v)=\int_0^1 \la\xi_1, v\ra \;\dd t =  \la\xi
_1, v\ra,
\]
which shows that $\lambda$ is the identity.

\subsection*{Step 2}
Assume that $(G\toto M, \omega +\phi )$ is an $\alpha$-simply
connected  non-degenerate twisted symplectic groupoid.
Let $\pi$ be its induced  $\phi$-twisted Poisson 
structure on  $M$. Then the above integration
process integrates the Lie algebroid $T^*M_{(\pi , \phi )}$
into a Lie groupoid, which is known to be isomorphic
to $G\toto M$, and a multiplicative $2$\ndash form $\omega'$
on that groupoid. By identifying this 
groupoid with $G\toto M$, therefore one may
think $\omega'$ as a multiplicative $2$\ndash form on $G$.
One needs to show that $\omega' =\omega$. By
Step 1, we conclude that $\omega' $ and $\omega$ must
coincide along the identity space $M$.
Let $\tomega =\omega -\omega'$. Then $\tomega$ is
a multiplicative closed $2$\ndash form on $G$ and $\tomega|_{M}=0$.
 Given any $\xi \in \gm (A)$,  it is easy to see that
$(\Vec{\xi}(\alpha (x) ), 0, \Vec{\xi} (x))$ is
tangent to the graph $\Lambda$ of groupoid multiplication.
On the other hand, for any $\delta_x \in T_x G$, it is
also clear that $(\alpha_* \delta_x , \delta_x , \delta_x)\in
T\Lambda$. It thus follows that
$$\tomega (\Vec{\xi}(\alpha (x) ) ,  \alpha_* \delta_x )-
\tomega ( \Vec{\xi} (x), \delta_x )=0. $$
Therefore we have $\Vec{\xi}\per \tomega =0$. Thus
$$L_{\Vec{\xi} }\tomega =(di_{\Vec{\xi} }+i_{\Vec{\xi}} d) \tomega=0,$$
which implies that $\tomega =0$ since any point in $G$ can be
reached by a product of (local)  bisections generated  by $\Vec{\xi}$.
This concludes the proof.




\end{document}